\documentclass[a4paper, 12pt]{article}

\usepackage[latin1]{inputenc}
\usepackage[T1]{fontenc}
\usepackage{amsmath}
\usepackage{amssymb}
\usepackage{amsthm}
\usepackage{paralist}

\normalsize
\providecommand{\GenBy}[1]{\bigl\langle #1 \bigr\rangle}

\title{A locally finite tree that behaves like an infinite star}
\author{Mykhaylo Tyomkyn}

\begin{document}

\maketitle

\begin{abstract}
We construct a tree $T$ of maximal degree~$3$ with infinitely many leaves such that whenever finitely many of them are removed, the remaining tree is isomorphic to $T$. In this sense $T$ resembles an infinite star. 
\end{abstract}

\section{Introduction and general remarks}
The famous Rado graph $R$ (see \cite{diestel},\cite{rado}) has the property of being isomorphic to $R-W$, where $W\subset V(R)$ is \emph{any} finite set of vertices. A similar property of a countable tree $T$ with infinitely many leaves would be $T \cong T-W$ for any finite set of \emph{leaves} $W$. Trivially, this property is satisfied by an infinite star, but what happens if we restrict ourselves to locally finite trees? If such a tree exists, it must be fairly ``large''. In this article we construct such a tree $T$ of maximal degree $3$, which is best possible.

The problem is motivated by the work of Bonato and Tardif~\cite{bonato} and the author~\cite{tyomkyn} on mutual embeddings of infinite trees. Call two trees \emph{twins} or \emph{mutually embeddable} if each of them contains an isomorphic copy of the other as a subtree. Define the \emph{twin number} of a tree $T$, in notation $m(T)$, to be the number of isomorphism classes of trees twinned with $T$. Bonato and Tardif~\cite{bonato} asked, what values the twin number can take, i.e. how many pairwise non-isomorphic trees can be mutually embeddable. Their conjecture was that the twin number of a tree is always either $1$ or infinite. In~\cite{tyomkyn} the author considered this problem for the class of locally finite trees. While it is easy to construct locally finite trees of twin number $1$, we could not find a locally finite tree $T$ of $m(T)=1$ which has a non-surjective self-embedding, other than the one-way infinite path. It can be shown that such a tree must be isomorphic to $T-x$ for infinitely many of its leaves $x$, so the natural question to ask is, if such a tree exists. The aim of this paper is to give an affirmative answer, and in fact the statement we prove is much stronger. Our construction might be helpful for solving the orginal problem and, we believe, it should be interesting in its own right.

The additional restriction of having infinitely many leaves is not substantial, as the only tree with finitely many leaves and the above property is a ray, i.e.\ a one-way infinite path. More dramatically, a ray is the only tree $T$ with finitely many leaves satisfying $T \cong T-x$ for \emph{some} leaf $x$. Indeed, the removal of a leaf from a tree with finitely many leaves either decreases the number of leaves by $1$ or does not affect it. In the latter case the vertex next to the removed leaf must have had degree $2$. Therefore, if $T$ has a vertex of degree at least $3$ and we start removing leaves of $T$ one by one in any order, we will at some point obtain a tree with fewer leaves than originally. In particular, repeated application of the isomorphism between $T$ and $T-x$ would result in a tree with fewer leaves than $T$ but isomorphic to $T$, a contradiction.

We now describe the construction of the desired tree $T$.  First of all, observe that it suffices to ensure $T \cong T-x$ for any single leaf $x$, for then the desired property follows by iteration. We build up $T$ inductively as a union of an ascending chain of trees $T_0 \subset T_1 \subset T_2 \subset \dotsb$ with $\Delta(T_n)= 3$ for all $n$. To be precise, we set $T = (V,E)$ where $V = \bigcup_{i=0}^{\infty}V(T_i)$, noting that the sequence $\{V(T_i)\}_{i=0}^{\infty}$ is ascending as well, and let $e \in E$ if and only if $e \in E(T_i)$ for some $i$. It is immediate that $T$ is another tree of maximal degree $3$.

Let the \emph{core} $c(T)$ of a tree $T$ be the graph spanned by all vertices $v$ such that $T-v$ has at least two infinite components. It is not hard to see that the core must be connected, i.e.\ $c(T)$ is a subtree of $T$. Note also, that the core is invariant under any isomorphism $\phi\colon T \rightarrow T-x$, i.e.\ $\phi$ restricted to $c(T)$ defines an automorphism of $c(T)$, which we also denote by $\phi$. Consequently, the inverse map $\phi^{-1}\colon T-x \rightarrow T$ gives rise to the inverse automorphism $\phi^{-1}\colon c(T) \rightarrow c(T)$.

Let $S$ and $T$ be two trees with vertices $v\in S$ and $w\in T$. The \emph{sum} of two rooted trees $(S,v)+(T,w)$ is constructed by identifying $v$ with $w$ and ``gluing together'' $S$ and $T$ accordingly. Similarly, if $W\subset V(T)$, we can define $(T,W)+(S,v)$ by attaching a copy of $S$ to each $w\in W$ identifying $v$ with $w$.

Let $x_1,x_2, \dotsc, x_n$ be leaves of a tree $T$ such that $T \cong T-x_i$ for each $i$ and let $\phi_i$ be the corresponding isomorphism between $T$ and $T-x_i$. Denote by $\Phi$ the set $\{\phi_1, \phi_2, \dotsc, \phi_n\}$. For $Y \subset c(T)$, define the \emph{closure} $\Phi(Y)$ of $Y$ under $\Phi$ as the set of all vertices $w \in c(T)$ that can be mapped to some $y \in Y$ by finitely applications of $\phi_1, \dotsc, \phi_n$ and $\phi_1^{-1}, \dotsc, \phi_n^{-1}$, that is, $\Phi(Y)$ is the image of $Y$ under $\GenBy{\Phi}$, the free group generated by $\Phi$. Furthermore, given a rooted tree $(S,s)$ we define $*(T,Y,S,s;\Phi)$, the \emph{convolution} of $(T,Y)$ and $(S,s)$ via $\Phi$, to be the tree $(T,\Phi(Y)) + (S,s)$.  As $Y \subset \Phi(Y)$, the convolution $*(T,Y,S,s;\Phi)$ extends the sum $(T,Y) + (S,s)$. More importantly, $*(T,Y,S,s;\Phi)-x_i \cong *(T,Y,S,s;\Phi)$ for all $i$ such that the underlying isomorphism extends the one between $T$ and $T-x_i$ and conversely, the convolution $*(T,Y,S,s;\Phi)$ is the minimal extension of $(T,Y) + (S,s)$ preserving the above isomorphism. 


The construction of the desired tree $T$ comprises two major parts. First we construct the sequence $\{T_n\}$. Then we show that its union $T$ has the needed property. 

\section{Constructing the sequence}

Our aim is to construct an ascending sequence of trees $T_0\subset T_1\subset T_2\subset \dotsb$ along with a sequence of vertices $x_0,x_1,x_2,\dotsc$ such that for all $n$ the following conditions are satisfied:
\begin{enumerate}
\item $\Delta(T_n)=3$.\label{itm:1}
\item $x_i$ is a leaf of $T_n$ for all $i \le n+1$.\label{itm:2}
\item $T_n \cong T_n - x_i$ for all $i \le n$. Moreover, for $i\le m<n$ the isomorphism between $T_n$ and $T_n-x_i$ extends the isomorphism between $T_m$ and $T_m-x_i$.\label{itm:3}
\item Writing $\phi^{n}_i$ for the isomorphism between $T_n$ and $T_n -x_i$, as in $3$. and $\Phi$ for the set $\{\phi^n_0,\dots,\phi^n_n\}$ , there exists a vertex $y \in c(T_n)$ such that each vertex $y' \in \Phi(y)$ has degree $2$ in $T_n$.\label{itm:4}
\item Among all leaves in $T_n$ apart from $x_0,x_1,\dotsc, x_n$ the leaf $x_{n+1}$ has the minimal distance from $x_0$.\label{itm:5}
\end{enumerate}

Let $T_0$ be a doubly infinite ray, labeled by the integers $\dotsc, -1,0,1,2, \dotsc$ with a leaf attached to each \emph{even} non-negative vertex.  If we let $x_0$ and $x_1$ be the leaves attached to the vertices labeled by $0$ and $2$, then a trivial check confirms properties~\ref{itm:1}--\ref{itm:5} for $T_0$.

Suppose that we have constructed $T_n$ along with $x_{n+1}$. We wish to extend $T_n$ to $T_{n+1}$. Note that, by property \ref{itm:3} and the fact that $T_n$ is not a ray (which is implied by property \ref{itm:1}) $T_n$ must have infinitely many leaves, as explained in the first section. Thus, once we have constructed $T_{n+1}$ such that properties \ref{itm:1}--\ref{itm:4} are satisfied, we can choose $x_{n+2}$ according to \ref{itm:5}. 

We want to construct $T_{n+1}$ as a union of another ascending sequence of trees $U_0 \subset U_1 \subset U_2 \subset \dotsb$ where $U_0 = T_n$ and for all $i$ the tree $U_i$ should have the following properties:
\begin{enumerate}[(i)]
\item $\Delta(U_i)=3$.\label{itm:i}
\item $x_0, x_1, \dots, x_{n+1}$ are leaves of $U_i$.\label{itm:ii}
\item If $i$ is even, then $U_i \cong U_i - x_j$ for $j = 0, 1, \dots n$. If $i$ is odd, then $U_i \cong U_i - x_{n+1}$. In both cases the corresponding isomorphism between $U_i$ and $U_i - x_j$ extends the one between $U_{i-2}$ and $U_{i-2}-x_j$.\label{itm:iii}
\end{enumerate}

At the same time we construct a sequence of vertex sets $W_i \subset U_{i-1}$, whose meaning will become clear later on.

Choose a vertex $y$ as in \ref{itm:4}. Attach a new path $P$ of length $2$ to $y$, whose other end we denote by $y'$. Attach a new doubly infinite path $R$ to $y'$. Let $R_1$ and $R_2$ denote the two rays, into which $y'$ divides $R$. Now to each vertex on $R_1$ except $y'$ attach a copy of $(P,y)+(T_n,y)$ in the same way as it is attached to $y'$. To each vertex on $R_2$ attach a copy of $(P,y)+(T_n-x_{n+1},y)$. Call the resulting tree $U_1$ and let $W_1 = \{y\}$ and $S = U_1-(U_0-y)$. In other words,  $U_1 = (U_0,y)+(S,y)$.

To obtain $U_{2i}$ from $U_{2i-1}$, set $U_{2i} = *(U_{2i-2},W_{2i-1},S,y;\Phi_{2i-2})$ where $\Phi_{2i-2}$ is the set of isomorphisms removing $x_0, \dotsc, x_n$ from $U_{2i-2}$ (note that $U_{2i-1} = (U_{2i-2},W_{2i-1})+(S,y)$ so $U_{2i}$ extends $U_{2i-1}$). Define $W_{2i}$ to be $\Phi_{2i-2}(W_{2i-1})\setminus W_{2i-1}$, in other words $W_{2i}$ is the set of those vertices in $U_{2i-1}$ to which we attached a new copy of $(S,y)$.

Similarly, to obtain $U_{2i+1}$ from $U_{2i}$, set $U_{2i+1}= *(U_{2i-1},W_{2i},S,y;\Phi_{2i-1})$ where $\Phi_{2i-1}$ consists of the single isomorphism between $U_{2i-1}$ and $U_{2i-1}-x_{n+1}$. Again, $U_{2i} = (U_{2i-1},W_{2i})+(S,y)$, and hence $U_{2i+1}$ extends $U_{2i}$. Define $W_{2i+1}$ to be $\Phi_{2i-1}(W_{2i})\setminus W_{2i}$.

Note that $W_1 \subset c(U_0)$ and therefore, by the invariance of the core, $W_2 \subset \Phi_0(W_1) \subset c(U_0) \subset c(U_1)$. Hence, $W_3 \subset \Phi_1(W_2) \subset c(U_1)\subset c(U_2)$. It follows by induction that $W_i \subset c(U_{i-1})$ for all $i$. This shows that a leaf of $U_i$ remains a leaf in $U_{i+1}$. In particular, $x_0, x_1, \dotsc, x_{n+1}$ are leaves of all $U_i$ and (\ref{itm:ii}) holds.

Another important observation is the fact that, for even $i$, the set $W_i$ lies on the same side of $P$ as $y$, whereas for odd $i \geq 3$, $W_i$ lies on the same side as $y'$.

We also note that property~(\ref{itm:iii}) is immediate from the construction and the corresponding property of the convolution.
 
We now must prove (\ref{itm:i}). For $U_1$ this follows by construction. Property~\ref{itm:4} of $T_n$ gives us $\Delta(U_2) = 3$. To prove the general statement it suffices to show that $d_{U_i}(w) = 2$ for all $w \in W_{i+1}$. Indeed, since we know that $d_S(y) = 1$, this will imply that at each step we identify vertices of degree~$2$ with vertices of degree~$1$, i.e.\ no vertex of degree~$4$ or greater is generated. If $w \in W_{i+1}$, then $w = \phi(w')$ for some $w' \in W_{i}$ and $\phi \in \GenBy{\Phi}$ (group generated by $\Phi$). Therefore, by the induction hypothesis, $d_{U_i}(w)=d_{U_i}(w')=2$ unless $w \in \Phi_{i-1}(v)$ for some $v$ adjacent to one of $x_0, \dots, x_n$ if $i$ is odd and adjacent to $x_{n+1}$ if $i$ is even.
Note, however, that $W_3 \not \subset U_0$, which means that for odd $i$ we have $W_i \not \subset U_0$, so $w$ or $w'$ cannot be mapped to a vertex adjacent to $x_0,\dotsc,x_n$. And if $i$ is even, the above obstruction can only occur in the case $i=2$, however, by construction of $U_1$ we know that this cannot happen. Therefore, property~(\ref{itm:i}) holds.


We need to show that $T_{n+1}$, defined as the union of $U_0 \subset U_1 \subset \dotsb$, satisfies \ref{itm:1}--\ref{itm:4}. Property~\ref{itm:1} is immediate from~(\ref{itm:i}). Property~\ref{itm:2} follows from (\ref{itm:ii}) and the fact that we can always choose $x_{n+2}$ according to \ref{itm:5}. Property \ref{itm:3} is a consequence of (\ref{itm:iii}).
So the only property left to be checked is \ref{itm:4}.

Let $z$ be the center of $P$, i.e.\ the vertex between $y$ and $y'$. Note that, by construction of $T_{n+1}$, vertices in $\Phi(z)$ cannot have degree~$1$ in $T_{n+1}$. If some $z' \in \Phi(z)$ has degree~$3$ then $z'$ has degree~$3$ in some $U_i$. But this implies that $z$ has degree~$3$ in some $U_j$, which contradicts the fact that $z \notin W_j$ for any $j$. 

\section{Taking the union}

We now define $T$ to be the union of the ascending sequence~$\{T_n\}$. We already mentioned that $\Delta(T)=3$. Note also that, by property~\ref{itm:3} of the sequence $\{T_n\}$, $x_i$ is a leaf of $T$ for all $i$ and property~\ref{itm:5} implies that $T$ has no other leaves. Finally, by property~\ref{itm:3}, $T$ is isomorphic to $T-x_i$ for all $i$. Indeed, since the isomorphisms $\phi_i^n\colon T_n \rightarrow T_n-x_i$ are ``nested'', we can combine them and define $\phi_i\colon T \rightarrow T-x_i$ by $\phi_i(y)=\phi_i^n(y)$, where $n$ is the smallest index satisfying $y \in T_n$. Hence, we have shown that $T$ is a tree with the desired properties.

\section{Further remarks}

The above construction seems to have a very high ``degree of freedom'', i.e.\ altering the construction one can obtain many pairwise nonisomorphic trees with the above property. It would be interesting to find out what properties a locally finite tree can have in addition to $T$ with $T \cong T-x$ for each leave $x$ . In particular, we do not know, whether such a tree $T$ must or can be isomorphic to $T-X$ for some \emph{infinite} set of leaves $X$. This is closely related to the problem of twin numbers and the solution could shed more light on ``paradoxical'' properties of graphs.

\end{document}